# The Z Transform over Finite Fields

R. M. Campello de Souza, H. M. de Oliveira, D. Silva[1]

*Abstract* — Finite field transforms have many applications and, in many cases, can be implemented with a low computational complexity. In this paper, the Z Transform over a finite field is introduced and some of its properties are presented.

## I. INTRODUCTION

Discrete transforms play a very important role in Engineering. Particularly significant examples are the well known Discrete Fourier Transform (DFT) and the Z Transform [1], which have found many applications in several areas, specially in the field of Electrical Engineering. A DFT over finite fields was also defined [2] and applied as a tool to perform discrete convolutions using integer arithmetic. Recently, a discrete Hartley transform over finite fields was introduced [3], which has interesting applications in the fields of digital multiplexing and spread spectrum [4], [5].

In this paper, the finite field Z transform is introduced. Earlier attempts to deal with Z transforms over finite fields considered only finite sequences of elements of a Galois field, a clever stratagem to evade the problem of series convergence over a finite algebraic structure [6]. In order to define a Z-transform over a finite field, we need to answer a few questions such as: *What is the meaning of a complex variable z defined over a finite field? How can we represent the Finite-Field Argand-Gauss plane? What about the convergence of infinite series over finite fields? Are there convergent series? In which sense?* Simple questions such as

$$\sum_{k=0}^{+\infty} 2^{-k} := \sum_{k=0}^{+\infty} [2^k]^{-1} (mod\ 7) \equiv 1+4+2+1+4+2+1+4+... = ?$$

$$\sum_{k=0}^{+\infty} 5^{-k} := \sum_{k=0}^{+\infty} [5^k]^{-1} (mod\ 7) \equiv 1+3+2+6+4+5+1+3+... = ?$$

should be addressed.

In the next section some mathematical preliminaries are presented. In particular, Gaussian integers over GF($p$) are defined and the unimodular group of GF($p^2$) is constructed. Section 3 introduces a representation for the complex Z plane over a Galois field GF($p$). Infinite sequences over finite fields are discussed in section 4, where the problem of series convergence is investigated. In section 5 the finite field Z transform is introduced and some of its properties are derived. The discrete-time finite field Fourier transform is also defined. The paper closes with a few concluding remarks.

## II. MATHEMATICAL PRELIMINARIES

II.1 Complex numbers over finite fields

The set G($p$) of Gaussian integers over GF($p$) defined below plays an important role in the ideas introduced in this paper (hereafter the symbol := denotes *equal by definition*). **Q**, **R** and **C** denote the rational, real and complex sets, respectively. $\delta$ is the Kronecker symbol.

**Definition 1**: G($p$) := {$a + jb$, $a, b \in$ GF($p$)}, $p$ being an odd prime for which $j^2 = -1$ (i.e., $p \equiv 3\ (mod\ 4)$) is a quadratic non-residue in GF($p$), is the set of Gaussian integers over GF($p$).

Let $\otimes$ denote Cartesian product. It can be shown, as indicated below, that the set G($p$) together with the operations $\oplus$ and $*$ defined below, is a field [7].

**Proposition 1**: *Let*

$\oplus$ : G($p$) $\otimes$ G($p$) $\to$ G($p$)
($a_1+jb_1, a_2+jb_2$) $\to$ ($a_1+jb_1$) $\oplus$ ($a_2+jb_2$) =
    =($a_1+a_2$)+$j(b_1+b_2)$,
*and*
$*$ : G($p$) $\otimes$ G($p$) $\to$ G($p$)
($a_1+jb_1, a_2+jb_2$) $\to$ ($a_1+jb_1$) $*$ ($a_2+jb_2$)=
    =($a_1a_2 - b_1b_2$) + $j(a_1b_2+a_2b_1)$.
*The structure* GI($p$) := < G($p$), $\oplus$, $*$ > *is a field.*
*In fact,* GI($p$) *is isomorphic to* GF($p^2$). □

By analogy with the complex numbers, the elements of GF($p$) and of GI($p$) are said to be real and complex, respectively.

**Definition 2**: *The unimodular set of* GI($p$), *denoted by* $G_1$, *is the set of elements* $\zeta=(a+jb) \in$ GI($p$), *such that* $a^2+b^2 \equiv 1\ (mod\ p)$.

**Proposition 2**: $\zeta^{p+1} \equiv |\zeta|^2 \equiv a^2 + b^2 \pmod{p}$.

**Proof**: $\zeta^p = (a+jb)^p \equiv a^p + j^p b^p \pmod{p}$, once GI($p$) is isomorphic to GF($p^2$), a field of characteristic $p$. Since $p=4k+3$, $j^p = -j$, so that $\zeta^p \equiv a - jb\ (mod\ p) = \zeta^* (mod\ p)$. Therefore, $\zeta^{p+1} \equiv \zeta\zeta^* = |\zeta|^2 \equiv a^2 + b^2 \pmod{p}$. □

**Proposition 3**: *The structure* <$G_1$, $*$> *is a cyclic group of order* ($p+1$).

**Proof**: $G_1$ is closed with respect to multiplication. On the other hand, it is a well know fact that the set of non-nzero elements of GF($p$), with the operation of multiplication, is a cyclic group of order ($p$-1) (denoted here by G) [8]. Therefore $G_1$ is a cyclic subgroup of G and, from proposition 2, it has order $p+1$. □

To determine the elements of the unimodular group it helps to observe that if $\zeta=a+jb$ is one such element, then so is every element in the set $\Gamma$={$b+ja$, ($p-a$)+$jb$, $b+j(p-a)$, $a+j(p-b)$, ($p-b$)+$ja$, ($p-a$)+$j(p-b)$, ($p-b$)+$j(p-a)$}.

**Example 1**: *Unimodular groups of* GF($7^2$) *and* GF($11^2$). In each case, table I lists the elements of the subgroups $G_1$ of order 8 and 12, and their orders.

TABLE I. ELEMENTS OF $G_1$.

| $\zeta \in$ GI(7) | Order | $\zeta \in$ GI(11) | Order |
|---|---|---|---|
| 1 | 1 | 1 | 1 |
| -1 | 2 | -1 | 2 |
| $j, -j$ | 4 | $5+j3, 5+j8$ | 3 |
| $2+j2, 2+j5, 5+j2, 5+j5$ | 8 | $j, -j$ | 4 |
| | | $6+j8, 6+j3$ | 6 |
| | | $8+j6, 8+j5, 3+j6, 3+j5$ | 12 |

Figure 1 illustrates the 12 roots of unity in GF($11^2$). Clearly, $G_1$ is isomorphic to $C_{12}$, the group of proper rotations of a

Departamento de Eletrônica e Sistemas - CTG - UFPE,
Communications Research Group, C.P. 7800, 50711 - 970, Recife - PE, Brasil
E-mail: {ricardo, hmo}@ufpe.br, danilos@hotlink.com.br

regular dodhecaedron. $\zeta = 8+j6$ is a group generator corresponding to a counter-clockwise rotation of $2\pi/12 = 30°$. Symbols of the same colour indicate elements of same order, which occur in conjugate pairs.

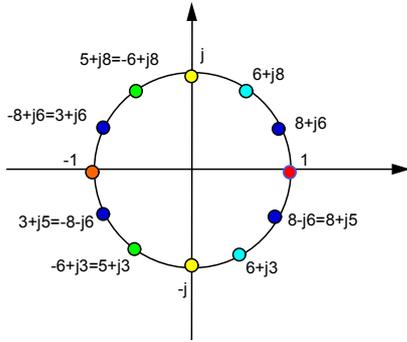

**Figure 1.** Roots of unity in $GF(11^2)$ expressed as elements of $GI(11)$.

II.2 *Polar Form for Gaussian Integers in a Finite Field*

It is a well known fact that, in the usual complex number arithmetic, the so-called polar representation has many interesting aspects that make it very attractive for many applications, particularly when operations such as multiplication and exponentiation are necessary. Keeping that in mind and aiming to create a representation that will make possible more efficient implementations of arithmetic modulo $p$, a new representation for the elements of $GI(p)$ is introduced in this section.

In the definition of $GI(p)$, the elements were written in rectangular form $\zeta = (a+jb)$. In what follows, a different representation for the elements of the multiplicative group of $GI(p)$ is proposed, which allows to write them in the form $r\varepsilon^\theta$. By analogy with the continuum, such a form is said to be polar.

**Proposition 4**: *Let $G_A$ and $G_B$ be subgroups, of the multiplicative group $G_C$ of the nonzero elements of $GI(p)$, of orders $N_A = (p-1)/2$ and $N_B = 2(p+1)$, respectively. Then all nonzero elements of $GI(p)$ can be written in the form $\zeta = AB$, where $A \in G_A$ and $B \in G_B$.*

**Proof**: Since $G_C$ is a cyclic group and both, $N_A$ and $N_B$, are divisors of $p^2-1$, then the subgroups $G_A$ and $G_B$ of $GI(p)$ do exist. Now, the direct product $(G_A \otimes G_B)$ [8] has order $p^2-1$, because, since $p$ is of the form $4k+3$, then the greatest common divisor between $N_A$ and $N_B$ satisfies $GCD(N_A,N_B) = GCD(2k+1, 4(2k+2)) = 1$; i.e., $|(G_A \otimes G_B)| = lcm(|G_A|,|G_B|) = N_A N_B = p^2-1$. Therefore, $(G_A \otimes G_B)$ is the multiplicative group of $GI(p)$, and every element of $GF(p)$ can be written in the form $\zeta = \alpha\beta$, $\alpha \in G_A$, $\beta \in G_B$. □

Considering that any element of a cyclic group can be written as an integral power of a group generator, it is possible to set $r = \alpha$ and $\varepsilon^\theta = \beta$, where $\varepsilon$ is a generator of $G_B$. The powers $\varepsilon^\theta$ of this element play, in some sense, the role of $e^{j\theta}$ over *the* complex field. Thus, the polar representation assumes the desired form, $\zeta = r\varepsilon^\theta$. Note that $2(p+1)$ plays the role of $2\pi$.

In the above discussion, it seems clear that $r$ is going to play the role of the modulus of $\zeta$. Therefore, before further exploring the polar notation, it is necessary to define formally the concept of modulus of an element in a finite field. Considering the nonzero elements of $GF(p)$, it is a well-known fact that half of them are quadratic residues of $p$ [9]. The other half, those that do not possess square root, are the quadratic nonresidues. Likewise, in the field **R** of real numbers, the elements are divided into positive and negative numbers, which are, respectively, those that possess and do not possess a square root. The standard modulus operation in **R** always gives a positive result. By analogy, the modulus operation in $GF(p)$ is going to be defined, such that it always results in a quadratic residue of $p$.

**Definition 3**: *The modulus of an element $a \in GF(p)$, where $p=4k+3$, is given by*

$$|a| := \begin{cases} a, & \text{if } a^{\frac{p-1}{2}} \equiv 1 \ (mod \ p) \\ -a, & \text{if } a^{\frac{p-1}{2}} \equiv -1 \ (mod \ p). \end{cases}$$

**Proposition 5**: *The modulus of an element of $GF(p)$ is a quadratic residue of $p$.*
**Proof**: Since $p=4k+3$, it implies that $(p-1)/2 = 2k+1$, such that $(-1)^{(p-1)/2} \equiv -1 \ (mod \ p)$. By Euler's criterion [9], if $a^{(p-1)/2} \equiv 1 \ (mod \ p)$, then $a$ is a quadratic residue of $p$; if $a^{(p-1)/2} \equiv -1 \ (mod \ p)$, then $a$ is a quadratic nonresidue of $p$. Therefore, $(-a)^{(p-1)/2} \equiv (-1)(-1) \equiv 1 \ (mod \ p)$, and it follows that $a$ is a quadratic residue of $p$. □

**Definition 4**: *The modulus of an element $a+bj \in GI(p)$, where $p=4k+3$, is given by $|a+jb| := \left| \sqrt{|a^2+b^2|} \right|$.*

The inner modulus sign in the above expression is necessary in order to allow the computation of the square root of the quadratic norm $a^2+b^2$, and the outer one guarantees that such an operation results in one value only. In the continuum, such expression reduces to the usual norm of a complex number, since both, $a^2+b^2$ and the square root operation, produce only positive numbers.

At this point it is convenient to substitute the groups $G_A$ and $G_B$ denominations for ones that are more appropriate to the polar representation.

**Definition 5**: *The group of modulus of $GI(p)$, denoted by $G_r$, is the subgroup of order $(p-1)/2$ of $GI(p)$.*

**Definition 6**: *The group of phases of $GI(p)$, denoted by $G_\theta$, is the subgroup of order $2(p+1)$ of $GI(p)$.*

**Proposition 6**: *If $\zeta = a+jb = r\varepsilon^\theta$, where $r \in G_r$ and $\varepsilon^\theta \in G_\theta$, then $r = |\zeta|$.*
**Proof**: Every element of $G_r$ has an order that divides $(p-1)/2$. Thus, if $r \in G_r$, then $r^{(p-1)/2} \equiv 1 \ (mod \ p)$, and $|r| = r$. Besides that, as shown in next section, the elements of the group $G_\theta$ are those $a+jb$ such that $a^2+b^2 \equiv \pm 1 \ (mod \ p)$. Therefore, according to definition 3, such elements have modulus equal to 1, which means that $|\zeta| = |r\varepsilon^\theta| = |r||\varepsilon^\theta| = r1 = r$. □

An expression for the phase $\theta$ as a function of $a$ and $b$ can be found by normalising the element $\zeta$ (that is, calculating $\zeta/r = \varepsilon^\theta$), and then solving the discrete logarithm problem of $\zeta/r$ in the base $\varepsilon$ over $GF(p)$. Thus, the conversion rectangular to polar form is possible. The inverse one is done simply by the exponentiation operation. From the above it can be observed that the polar representation being introduced is consistent with the usual polar form defined over the complex numbers. The modulus belongs to $GF(p)$ (the modulus is a real number) and is a quadratic residue (a positive number), and the exponential component $\varepsilon^\theta$ has modulus one and belongs to $GI(p)$ ($e^{j\theta}$ also has modulus one and belongs to the complex field).

## III. THE Z PLANE

To define a Z transform over a finite field, it is necessary first to establish what is meant by the complex Z plane in GF($p$). In order to do so, it is necessary to introduce a special family of finite groups.

**Definition 7**: *The supra-unimodular set of* GI($p$), *denoted* $G_S$, *is the set of elements* $\zeta=(a+jb)\in$ GI($p$), *such that* $(a^2+b^2)^2\equiv 1$ (*mod p*).

**Proposition 7**: *If* $\zeta=a+jb$, *then* $\zeta^{2(p+1)}\equiv(a^2+b^2)^2$ (*mod p*).

**Proof**: $\zeta^p=(a+jb)^p\equiv a^p+j^p b^p$ (*mod p*), since GI($p$) is isomorphic to GF($p^2$), a field of characteristic $p$. Also, since $p=4k+3$, $j^p\equiv -j$, so that $\zeta^p\equiv a-jb$ (*mod p*), which means that $\zeta^{p+1}\equiv(a+jb)(a-jb)$ (*mod p*). Therefore, $\zeta^{2(p+1)}\equiv(a^2+b^2)^2$ (*mod p*). □

**Proposition 8**: *The structure* <$G_S$,*>, *is a cyclic group of order* $2(p+1)$, *called the supra-unimodular group of* GI($p$).

**Proof**: $G_S$ is closed with respect to *, since that if $(a+bj)$ and $(c+jd)$ belong to $G_S$, i.e., if $(a^2+b^2)^2\equiv(c^2+d^2)^2\equiv 1$ (*mod p*), then $e+jf=(a+jb)*(c+jd)=(ac-bd)+j(ad+bc)$. So that $(e^2+f^2)^2=(a^2c^2-2abcd+b^2d^2+a^2d^2+2abcd+b^2c^2)^2=((a^2+b^2)*(c^2+d^2))^2=(a^2+b^2)^2*(c^2+d^2)^2\equiv 1$ (*mod p*), which means that $(e+jf)\in G_S$. Now, $G_S$ is a closed subset of a cyclic group (the multiplicative group of GI($p$)), therefore $G_S$ is a cyclic subgroup. Besides that, from proposition 2, $\zeta\in G_S$ satisfies $\zeta^{2(p+1)}\equiv 1$ (*mod p*). Thus, $\zeta$ is one of the $2(p+1)$th roots of unity in GI($p$). There exists $2(p+1)$ such roots e therefore $G_S$ has order $2(p+1)$). □

The elements $\zeta=a+jb$ of the supra-unimodular group $G_S$ satisfy $(a^2+b^2)^2\equiv 1$ (*mod p*), i.e., $a^2+b^2\equiv \pm 1$ (*mod p*), and all have modulus 1, as in the unimodular group $G_1$. However, $G_S$ has a larger order than $G_1$. It is important to observe that, due to the fact that a cyclic group has only one subgroup of a given order [8], $G_S$ is precisely the group of phases $G_\theta$ introduced in definition 6. The problem of finding a group generator for $G_S$, is dealt with in proposition 9 below.

**Proposition 9**: *If p is a Mersenne prime* ($p=2^n-1$, $n>2$), *the elements* $\zeta=a+jb$ *such that* $a^2+b^2\equiv -1$ (*mod p*) *are the generators of* $G_S$.

**Proof**: Let $\zeta\in G_S$ have order $N$. Since $a^2+b^2\equiv -1$ (*mod p*), $N$ divides $2(p+1)=2^{n+1}$. However, $\zeta$ is not unimodular, so that $N$ does not divide $p+1=2^n$. Therefore, $N=2^{n+1}=2(p+1)$, and $\zeta$ is a generator of $G_S$ □

**Example 2**: Let $p=31$, a Mersenne prime, and $\zeta=6+j16$. From definition 3: $r=\left|\sqrt{6^2+16^2}\right|\equiv\left|\sqrt{13}\right|\equiv 7$ (*mod 31*),

so that $\varepsilon = \zeta/r=23+j20$ and $a^2+b^2=23^2+20^2\equiv -1$ (*mod* 31). Therefore $\varepsilon$ has order $2(p+1)=64$ (a generator). A unimodular element $\beta$ of order $N$, such that $N \mid 2^5$, can be found taking $\beta=\varepsilon^{2(p+1)/N}=\varepsilon^{64/N}$. o

A generator $\varepsilon$ of the supra-unimodular must be used to construct the Z plane over GF($p$). From the powers of $\varepsilon$ those elements that are on the unit circle are obtained. By multiplying those by the members of the group of modules, the remaining elements of GI(7) are identified on other circles on the Z plane.

**Proposition 10**: *In the Z complex plane over* GF($p$) *there are* $2(p+1)$ *elements in each circle*.

**Proof**: There are as many circles in the Z plane as the order of the group of modulus $G_r$. Therefore, denoting by $m$ the number of elements of a given modulus, it is possible to write $m|G_r|=p^2-1$ and the result follows. □

**Example 3**: The Z plane over GF(7). Let $p=7$, and $\zeta=6+j4$. From definition 3,

$$r=\left|\sqrt{6^2+4^2}\right|\equiv\left|\sqrt{3}\right|\equiv\left|\sqrt{4}\right|\equiv|2|\equiv 2 \ (mod\ 7),$$

so that $\varepsilon=\zeta/r=3+j2$ and $a^2+b^2=13\equiv -1$ (*mod* 31). Therefore $\varepsilon$ has order $2(p+1)=16$, so it is a generator of the group $G_S$. The Z plane over GF(7) is depicted in figure 2 below. The nonzero elements of GF($p$), namely $\pm 1, \pm 2, \pm 3$, are located on the horizontal axis, in the right or left side, according if they are, respectively, quadratic residues or quadratic non-residues of $p=7$ (respectively, positive or negative numbers, in the usual complex Z plane). There are three circles, of radius 1,2 and 4, corresponding to the $(p-1)/2=3$ elements of the group of modules $G_r$. A similar situation occurs for the elements of GI(7) of the form $jb$ (corresponding to imaginary elements in the continuous case). The 16 elements on the unit circle correspond to the elements of $G_S$ and are obtained as powers of $\varepsilon$. The even powers correspond to the elements of $G_1$ ($a^2+b^2\equiv 1$ (*mod* 7)) and the odd powers to the elements satisfying $a^2+b^2\equiv -1$ (*mod* 7). The remaining 32 elements of the Z plane are obtained simply by multiplying those on the unit circle by the modulus 2 and 4. Observe that the elements on the straight line $y=\pm x$ over a finite field also possesses the usual interpretation associated to $tg\theta=\pm 1$. Table II relates the points of figure 2 with their respective orders as elements of GI(7).

TABLE II.
THE NUMBER OF ELEMENTS OF A GIVEN ORDER IN THE Z PLANE OVER GF(7).

| element | ○ | ● | ● | ● | ● | ● | ● | ● | ● | ● |
|---|---|---|---|---|---|---|---|---|---|---|
| order | 1 | 2 | 3 | 4 | 6 | 8 | 12 | 16 | 24 | 48 |
| number of elements | 1 | 1 | 2 | 2 | 2 | 4 | 4 | 8 | 8 | 16 |

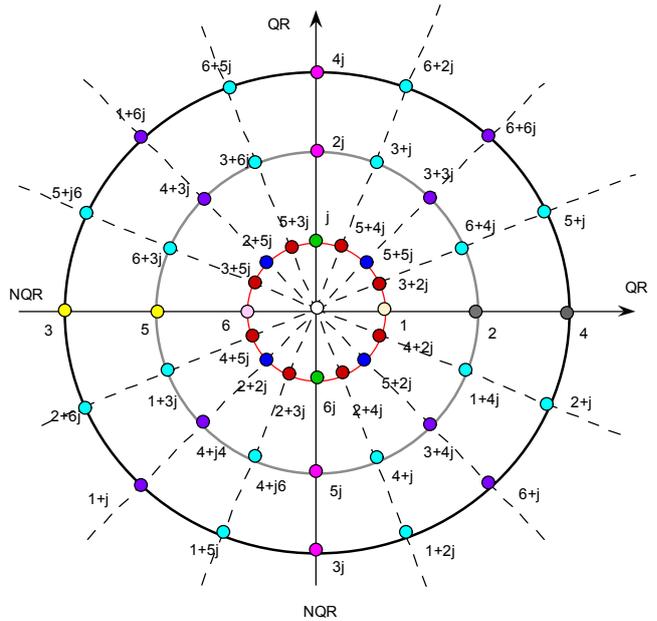

**Figure 2**. The Z Plane over the Galois Field GF(7).

When calculating the order of a given element, the intermediate results generate a trajectory on the Z plane, called the order trajectory. In particular, If ζ has order *N*, the trajectory goes through *N* distinct points on the Z plane, moving in a pattern that depends on *N*. Specifically, the order trajectory touches on every circle of the Z plane (there are $\|G_r\|$ of them), in order of crescent modulus, always returning to the unit circle. If it starts on a given radius, say *R*, it will visit, counter-clockwise, every radius of the form $R+kr$, where $r=(p^2-1)/N$ and $k=0,1,2,....,N-1$.

**Example 4**: Table III lists some elements ζ ∈ GI(7) and their orders *N*. Figures 3-5 show the order trajectories generated by ζ.

TABLE III.
SOME ELEMENTS AND THEIR ORDERS IN GI(7).

| ζ | 2*j* | 3+3*j* | 6+4*j* |
|---|---|---|---|
| *N* | 12 | 24 | 48 |

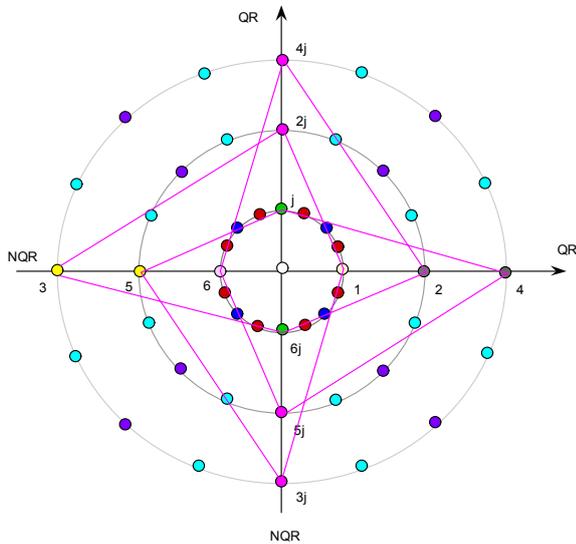

**Figure 3**. Order trajectory for ζ=*j*2, an element of order *N*=12 of GI(7), on the Z Plane over GF(7).

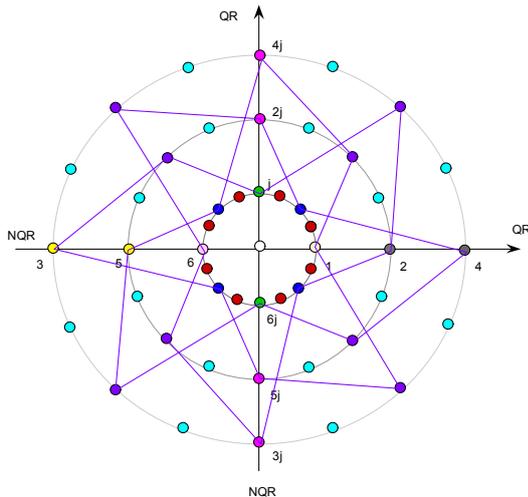

**Figure 4**. Order trajectory for ζ=3+*j*3, an element of order *N*=24 of GI(7), on the Z Plane over GF(7).

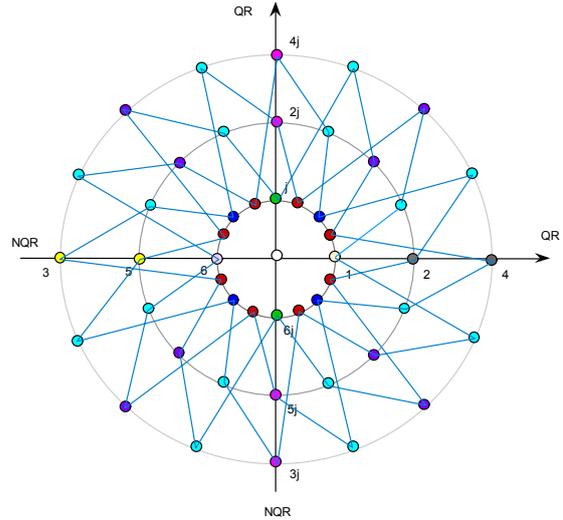

**Figure 5**. Order trajectory for ζ=6+*j*4, an element of order *N*=48 of GI(7), on the Z Plane over GF(7).

## IV. INFINITE SERIES OVER FINITE FIELDS

It is straightforward to conceive infinite sequences whose elements belong to a finite field. For instance, given a sequence of integers $\{x[n]\}_{-\infty}^{\infty}$, it is possible to generate a sequence over GF(*p*) simply by considering $\{x[n] \ (mod \ p)\}_{-\infty}^{\infty}$. In general, *x[n]* may even be a sequence of rational elements: any element $r/s \in \mathbf{Q}$ can be mapped over GF(*p*), assuming values $[r \ (mod \ p)].[s \ (mod \ p)]^{-1}$. The focus here concerns *periodic* infinite sequences over a finite field. For instance, consider the following right-sided GF(7)-valued sequences:

**Example 5**: Let $\{x[n]\}_{-\infty}^{\infty}$, $x[n] \in$ GF(*p*). For instance:

*i*) $\{3^n\}_{n=0}^{\infty}$ over GF(7), i.e., { 1 3 2 6 4 5 1 3 2 6 4 5 ... }

*ii*) $\{1^n\}_{n=0}^{\infty}$ over GF(7), i.e., { 1 1 1 1 1 1 1 1 1 1 1 1 ... }

We wonder whether a Z-transform of this sequence can be defined by $X(z) := \overset{?}{\underset{n=-\infty}{\overset{+\infty}{\sum}}} x[n]z^{-n}$, $z \in$ GI(*p*). Is this series convergent? The main interest here is in the convergence of infinite series, because important transforms, such as Fourier and Z transforms, involve a sum of infinitely many terms. It is a well known fact that the infinite series

$\sum_{n=1}^{\infty}(-1)^{n+1}$ =1-1+1 -1 +1 ... diverges in the classical sense. However, Euler and others noticed that the arithmetic mean of the partial sums converges to 1/2. The partial sums of this series are $S_1$=1, $S_2$=0, $S_3$=1, $S_4$=0 ... and the arithmetic mean $\sigma_n := \frac{1}{n}\sum_{k=1}^{n}S_k$ forms a sequence ($\sigma_n$) that converges to 1/2. When a series converges in the sense that the arithmetic mean of the partial sums converges, it is said to be Cesàro-summable (Ernesto Cesàro (1859-1906)). Every convergent series in the usual sense is Cesàro-summable and the series sum (in the usual sense) is equal to the limit of the sequence of the partial sums arithmetic mean. That shows the Cesàro summability concept is useful, since it can make divergent series summable.

A (new) convergence criterion, suitable for series over finite fields, which is derived from the Cesàro summability, is now introduced. Given $\{x[n]\}_{1}^{\infty}$, the partial sums *S[n]* are defined according to: $S[n] := \sum_{k=1}^{n}x[k]$.

**Definition 8**: *The Cesàro sum over a finite field is defined by* $\sigma_n := \left(\frac{1}{n}\sum_{k=1}^{n} S[k]\right)$, *where $S[k] \in GF(p)$ are interpreted as integers.*  ○

If $\{x[n]\}_1^\infty$ is a periodic sequence over GF(p), so is $\{S[n]\}_1^\infty$. Let $P$ denote the period of the latter sequence. Therefore

$$\lim_{n\to\infty}\left(\frac{1}{n}\sum_{k=1}^{n}S[k]\right) = \lim_{n\to\infty}\left(\frac{1}{n}\sum_{k=1}^{\lfloor n/P \rfloor P}S[k]\right) + \left(\frac{1}{n}\sum_{k=n-n(\bmod P)}^{n}S[k]\right).$$

The second term vanishes and

$$\lim_{n\to\infty}\sigma_n = \lim_{n\to\infty}\frac{\lfloor n/P \rfloor}{n}\left(\sum_{k=1}^{P}S[k]\right).$$

But $\lim_{n\to\infty}\frac{\lfloor n/P \rfloor}{n} = \frac{1}{P}$

so that $\lim_{n\to\infty}\left(\frac{1}{n}\sum_{k=1}^{n}S[k]\right) = \frac{1}{P}\sum_{k=1}^{P}S[k]$ and finally

$$\left(\lim_{n\to\infty}\sigma_n\right)(\bmod\ p) \equiv \frac{1}{P(\bmod\ p)}\sum_{k=1}^{P}S[k]\ (\bmod\ p).$$

The heart of the matter is to take first the limit $n\to\infty$, and then evaluate the result after reducing it modulo $p$.

**Definition 9**: *A series over a finite field is said to be Cesàro convergent to $\sigma$ if and only if*

$$\sigma := \left(\lim_{n\to\infty}\sigma_n\right)(\bmod\ p) \in GF(p).$$

**Corollary**: *Every periodic series over a finite field, with a nonzero period, that is, $P\neq 0\ (\bmod\ p)$, is Cesàro convergent.*

**Example 5 (Revisited)**:
i) $\{S[k]\} = \{1\ 4\ 6\ 5\ 2\ 0\ 1\ 4\ 6\ 5\ 2\ 0\ 1\ 4\ ...\}$, $P=6\ (\bmod\ 7)$.
Therefore, $\sigma = \frac{1}{6}(1+4+6+5+2+0) \equiv 3\ (\bmod\ 7)$. Let us now return to the series $\sum_{k=0}^{+\infty}5^{-k} := \sum_{k=0}^{+\infty}[5^k]^{-1}(\bmod\ 7) =$

$\sum_{k=0}^{\infty}[5^{-1}]^k\ (\bmod\ 7) = \sum_{k=0}^{\infty}3^k\ (\bmod\ 7) \equiv 1+3+2+6+4+5+1+...$

If we write $\sum_{k=0}^{+\infty}Z^k \equiv \frac{1}{1-Z}$ over GF(7) and assume $Z=3$, we find

$\sum_{k=0}^{+\infty}3^k\ (\bmod\ 7) \stackrel{?}{\equiv} \frac{1}{1-3} \equiv (-2)^{-1} \equiv 3\ (\bmod\ 7)$.

This series over GF(7) converges exactly to the Cesàro sum $\sigma$!

ii) The series $\sum_{k=1}^{\infty}x[n]$ derived from $\{x[n]\}_0^{+\infty} = \{1\}_0^\infty$ is not convergent over GF(7). Remark that
$\{S[k]\} = \{1\ 2\ 3\ 4\ 5\ 6\ 0\ 1\ 2\ 3\ 4\ 5\ 6\ 0\ ...\}$, $P\equiv 0\ (\bmod\ 7)$ and $\sigma$ is undefined (unbounded). In this case, $\sum_{k=0}^{+\infty}Z^k \stackrel{?}{=} \frac{1}{1-Z}$ has $Z=1$ as a pole, so that $\sum_{k=0}^{+\infty}1^k$ diverges over GF(7). Now we are able to investigate the following: "*Given a periodic sequence $\{x[n]\}_{-\infty}^{+\infty}$ over GF(p), what is the region of Cesàro-convergence on the Z-Finite-Field plane for the series* $\sum_{n=-\infty}^{+\infty}x[n]Z^{-n}$, $Z\in GF(p^2)$?"

## V. THE FINITE FIELD Z-TRANSFORM

### V.1. Basic Sequences

The Z transform introduced in this paper deals with sequences $x[n]$, $-\infty<n<\infty$, defined over the Galois field GF(p), which are obtained from basic types of sequences: $\delta[n]$, $u[n]$, $A(a^n)$.
i) The finite field impulse (Galois impulse), denoted $\delta[n]$, is the sequence $x[n]$ defined by $x[n] = \delta[n] := \begin{cases} 1, & n=0 \\ 0, & n\neq 0.\end{cases}$

As it happens with real sequences, any finite field sequence $x[n]$ can be expressed as a sum of shifted and scaled Galois impulses.
ii) The unit step over GF(p) (Galois unit step) is given by

$$x[n] = u[n] := \begin{cases} \equiv 1\ (\bmod\ p), & n\geq 0 \\ 0, & n<0.\end{cases}$$

iii) The exponential sequence is $x[n]=A(a)^n$, $A$ and $a \in GF(p)$. This sequence is periodic with period $P$, which is the multiplicative order of $a\ (\bmod\ p)$.

### V.2. The Z transform

**Definition 10**: *The finite field Z transform (FFZT) of a sequence $x[n]$ over GF(p), is the GI(p)-valued function*

$$X(Z) := \sum_{n=-\infty}^{\infty} x[n]Z^{-n},\quad Z\in GI(p).$$

In the infinite series of the above definition, convergence is considered in the sense of definition 9. For any given sequence, the set of values of $Z\in GI(p)$ for which $X(Z)$ converges is called the region of convergence (ROC).

A class of important and useful finite field Z transform is that for which $X(Z)$ is a rational function $P(Z)/Q(Z)$, where $P(Z)$ and $Q(Z)$ are polynomials in $Z$. The roots of $P(Z)$ and $Q(Z)$ are called the zeros and poles of $X(Z)$, respectively. For the Z plane over GF(p), the ROC of $X(Z)$ corresponds to the elements of GI(p) that are not poles of $X(Z)$.

**Example 6**: *Right-sided exponential sequence over GF(p)*. Let $x[n]=a^n u[n]$, $a \in GF(p)$ and $u[n]$ being the Galois unit step. In this case, since $x[n]$ is non-zero only for $n \geq 0$,

$$X(Z) = \sum_{n=-\infty}^{\infty} a^n u[n]Z^{-n} = \sum_{n=0}^{\infty}(aZ^{-1})^n.$$

Computing the partial sums:
$$S_1 = 1$$
$$S_2 = 1 + aZ^{-1}$$
$$S_3 = 1 + aZ^{-1} + (aZ^{-1})^2$$
..................................................
$$S_{N-1} = 1 + aZ^{-1} + (aZ^{-1})^2 + (aZ^{-1})^3 + ........... + (aZ^{-1})^{N-2}$$
$$S_N = 1 + aZ^{-1} + (aZ^{-1})^2 + (aZ^{-1})^3 + .......... + (aZ^{-1})^{N-2} + (aZ^{-1})^{N-1}.$$

Now, denoting by $N$ the multiplicative order of $(aZ^{-1})$, the period of the sequence $S[k]$ is precisely $N$. Therefore, it is possible to write $\sigma_N = \frac{1}{N}\sum_{i=0}^{N-1}(N-i)(aZ^{-1})^i$, which is equal to

$$\sigma_N = \frac{(aZ^{-1})^N - 1}{aZ^{-1} - 1} - \frac{1}{N}\sum_{i=0}^{N-1}i(aZ^{-1})^i.$$ Since $(aZ^{-1})$ has multiplicative order $N$, $\sigma_N = -\frac{Z}{N}\sum_{i=0}^{N-1}ia^iZ^{-i-1}$, which is the same as $\sigma_N = \frac{Z}{N}\sum_{i=0}^{N-1}\frac{d}{dZ}(a^iZ^{-i})$, or

$$\sigma_N = \frac{Z}{N}\frac{d}{dZ}\left(\sum_{i=0}^{N-1}(aZ^{-1})^i\right).$$

After some direct manipulation, this is seem to be equal to $\sigma_N = \frac{1}{1-aZ^{-1}}$, which is the desired Z transform $X(Z)$. The only pole of $X(Z)$ is $Z=a$, so that the ROC is $Z \neq a$. The ROC in the finite field case has a different structure from the ROC of the usual Z-transform and, besides that, there exists many regions of convergence for a given analytical expression for $X(Z)$.

Evaluating the FFZT on $G_s$ results in the discrete-time Fourier transform over GF(p).

**Definition 11**: *The finite field discrete-time Fourier transform of a sequence x[n] over* GF(*p*) *is the* GI(*p*)-*valued function*

$$X(\varepsilon^\theta) := \sum_{n=-\infty}^{+\infty} x[n]\varepsilon^{-n\theta},$$

where $\varepsilon \in G_s$ has multiplicative order $2(p+1)$.

The above finite field Fourier transform is Cesàro-convergent if the ROC of the FFZT $X(Z)$ includes $G_s$. It should be compared with the previous definition by Pollard [2].

### V.3. The inverse FFZT

**Lemma 1**: $\sum_{Z \in GI(p)} Z^i = (p^2-1)\delta_{0,i}$.

**Proof**: For $i=1$ the summation is clearly equal to $p^2-1$. From proposition 4

$$\sum_{Z \in GI(p)} Z^i = \left(\sum_{r \in G_r}(r)^i\right)\left(\sum_{\theta \in G_\theta}(\varepsilon^\theta)^i\right).$$

But $\sum_{\theta \in G_\theta} \varepsilon^{\theta i} = \begin{cases} \frac{1-\varepsilon^{i2(p+1)}}{1-\varepsilon^i} \equiv 0 & i \neq 0 \\ 2(p+1) & i = 0 \end{cases}$, and therefore

$$\sum_{Z \in GI(p)} Z^i = ||G_r||.||G_\theta||\delta_{0,i}. \quad \square$$

**Theorem 1**: *The inverse finite field Z transform (IFFZT) is*

$$x[n] = \frac{1}{p^2-1} \sum_{Z \in GI(p)} X(Z)Z^n.$$

**Proof**: By definition, $X(Z) = \sum_{k=-\infty}^{\infty} x[k]Z^{-k}$. Multiplying both sides by $Z^n$ and summing over $Z \in GI(p)$, one obtains

$$\sum_{Z \in GI(p)} X(Z)Z^n = \sum_{Z \in GI(p)} \left(\sum_{k=-\infty}^{\infty} x[k]Z^{-k}\right)Z^n.$$

Interchanging the order of summation, the right-hand side (RHS) becomes

$$RHS = \sum_{k=-\infty}^{\infty} x[k]\left(\sum_{Z \in GI(p)} Z^{n-k}\right),$$

but the inner summation, from lemma 1, is nonzero only for $k=n$, so that

$$RHS = (p^2-1)x[n]$$

and the result follows. $\quad \square$

It is interesting to observe that, although the direct transform is an infinite summation, defined so to deal with infinite sequences, the IFFZT requires only a finite sum over the field GI(*p*). The Z transform introduced in this paper satisfies most properties of the usual Z transform defined over the complex field, such as linearity, time-shift, multiplication by exponential and so on.

## VI. CONCLUSIONS

This paper discusses the problem of constructing a Finite Field Z Transform, capable of handling infinite sequences over the Galois Field GF(*p*). To deal with such sequences, the concept of Cesàro convergence was used and it was shown that periodic sequences over GF(*p*) converge to the arithmetic mean of the Cesàro sums of the sequence. The Z plane over GF(*p*) was constructed and a Finite Field Z Transform (FFZT) was defined. A general inversion formula for the FFZT was determined. The finite field discrete-time Fourier transform was also introduced, which corresponds to the FFZT evaluated on the supramodular group (the unit circle). The transforms introduced here are also able of handling finite sequences and therefore they are a useful tool to deal with both FIR and IIR filters defined over finite algebraical structures.